*Dedicated to the bright memory of*
*Vitalii Yakovych Skorobohat'ko*

# RELATIVISTIC MECHANICS OF CONSTANT CURVATURE

**R. Ya. Matsyuk**  UDC 517.972+531.12: 530.12

We consider an inverse variational problem for the lines of constant curvature in (pseudo-)Euclidean two-, three-, and four-dimensional spaces. The accumulated results are physically meaningful in the case of relativistic mechanics of particles.

**Introduction**

The present paper is a survey of some achievements in the application of one conjecture of Vitalii Yakovych Skorobohat'ko concerning the construction of $n$-point geometries in Euclidean spaces (of different dimensions) and in the Minkowski space of the special theory of relativity. This conjecture was advanced in [3]. In the simplest interpretation, it can be formulated as the necessity of finding the principle of least action for the lines $x^\alpha(s)$ in the space, which have a constant first Frénet curvature $k$, i.e., satisfy the following equation of the third order $\frac{dk}{ds} = 0$. In any space, it is possible to consider both parametrized and nonparametrized curves. If a curve is given by functions $x^\alpha(\zeta)$, then, in the case where nonparametrized curves are the solutions of the differential equation, the equation obeyed by the vector function $x^\alpha(\zeta)$ must be parametrically invariant, i.e., must contain all solutions that differ by an arbitrary change of the parameter $\zeta$ along the curve. In order to deduce a differential equation of variational type (i.e., a vector equation for which one can, in principle, indicated one or another locally defined Lagrange function), it is reasonable to combine two criteria: a criterion of variationality and a criterion of symmetry of the equation under a group of space transformations. The results obtained in what follows have the same form for the Euclidean and pseudo-Euclidean spaces and the difference in formulas is observed only in the definition of the scalar product (i.e., in the signature of the metric tensor).

**1. Variationality and Parametric Invariance**

The criterion of variationality of a differential equation is reduced to certain restrictions imposed on the form of equation. Its detailed description can be found in [5]. In order that the set of expressions

Pidstryhach Institute for Applied Problems in Mechanics and Mathematics, National Academy of Sciences of Ukraine, Lviv, Ukraine; Lepage Research Institute, Prešov, Slovakia; e-mail: romko.b.m@gmail.com.

$$\mathcal{E}_\alpha\left(x^\beta, \frac{dx^\beta}{d\zeta}, \frac{d^2x^\beta}{d\zeta^2}, \frac{d^3x^\beta}{d\zeta^3}\right)$$

can be the left-hand side of a certain autonomous vector Euler–Poisson equation $\mathcal{E}_\alpha = 0$, it must have the form

$$\mathcal{E}_\alpha = \mathcal{A}_{\alpha\beta}\dddot{x}^\beta + \left(\ddot{x}^\gamma \partial_{\dot{x}^\gamma}\right)\mathcal{A}_{\alpha\beta}\ddot{x}^\beta + \mathcal{B}_{\alpha\beta}\ddot{x}^\beta + \mathcal{C}_\alpha, \qquad (1)$$

where the matrix $\mathcal{A}_{\alpha\beta}$ must be antisymmetric and both the matrices $\mathcal{A}_{\alpha\beta}$ and $\mathcal{B}_{\alpha\beta}$ and the column $\mathcal{C}_\alpha$ may depend on the variables $x^\beta$ and their derivatives $\dot{x}^\beta$.

If Eq. (1) is parametrically invariant, then it is possible to represent each its solution $x^i(t)$ as the plot of a certain vector function whose dimension is lower by one; in this case, the index $i$ runs through a set of values whose number is lower than the index $\alpha$ by one. If the space is $n$-dimensional, then it is possible to assume that the index $\alpha$ runs through the values $0, 1, \ldots, n-1$, whereas the index $i$ runs through the values $1, \ldots, n-1$. In this case, we redenote the variable $x^0$ by $t$. The functions $x^i(t)$ now satisfy the vector variational equation $E_i = 0$ whose left-hand side is similar to (1):

$$E_i = A_{ij}x'''^j + \left(x''^\ell \partial_{x'^\ell}\right)A_{ij}x''^j + B_{ij}x''^j + c_i. \qquad (2)$$

In this equation, the derivatives with respect to the variable $t$ are denoted by primes. Moreover, the antisymmetric matrix $A_{ij}$, the matrix $B_{ij}$, and the column $c_i$ depend not only on the variables $x^i$ and their derivatives $x'^i$ but also on the variable $t$. In both expressions (1) and (2), these variable coefficients must satisfy a system of partial differential equations that can be used to find these coefficients by adding (if necessary) the symmetry condition imposed on the variational equation. The symmetry condition reflects the invariance of the equation under the action of a group of (pseudo-)Euclidean space transformations in a sense that these transformations preserve the set of solutions of the equation. Note that, in this case, we do not speak about the invariance of the left-hand side of differential equation; actually, we state that the transformed left-hand side of the equation is expressed via left-hand side of the original equation. This condition, in turn, also gives a certain ("determining") system of partial differential equations for the variable coefficients of the required equation. Both conditions (of variationality and invariance) imposed either separately or together may have no solutions but, on the other hand, may correspond to a unique solution.

The variable coefficients in expressions (1) and (2) satisfy the following relations:

$$\mathcal{A}_{ij} = \frac{1}{t^2}A_{ij}, \quad \mathcal{B}_{ij} = \frac{1}{t}B_{ij}, \quad \mathcal{C}_i = \dot{t}c_i,$$

$$\mathcal{A}_{0i} = \frac{1}{t^2}A_{ij}x'^j, \quad \mathcal{B}_{0i} = -\frac{1}{t}B_{ij}x'^j,$$

$$\mathcal{B}_{00} = \frac{1}{t}B_{0i}x'^i, \quad \mathcal{C}_0 = -\dot{t}c_ix'^i,$$

$$\mathcal{E}_i = \dot{t}^3 \left[ A_{ij} x'''^j + \left( x''^j \partial_{x'^j} \right) A_{ij} x''^j + \frac{1}{\dot{t}} B_{ij} x''^j + \frac{1}{\dot{t}^3} c_i \right].$$

The expressions $\mathcal{E}_\alpha$ satisfy the Weierstrass condition:

$$\dot{t}\mathcal{E}_0 + \dot{x}^i \mathcal{E}_i = 0. \tag{3}$$

## 2. Two-Dimensional Space

In the two-dimensional space, it is impossible to construct a parametric-invariant variational equation because every antisymmetric matrix of size one in Eq. (2) is trivial (equal to zero). However, it is possible to get an equation of the form (1) whose solutions are geodesic lines of constant curvature [4]:

$$\frac{\varepsilon_{\alpha\beta} \dddot{x}^\beta}{(\dot{x}_\beta \dot{x}^\beta)^{3/2}} - 3 \frac{\dot{x}_\beta \ddot{x}^\beta}{(\dot{x}_\beta \dot{x}^\beta)^{5/2}} \varepsilon_{\alpha\beta} \ddot{x}^\beta + m \frac{(\dot{x}_\beta \dot{x}^\beta)\ddot{x}_\alpha - (\dot{x}_\beta \ddot{x}^\beta)\dot{x}_\alpha}{(\dot{x}_\beta \dot{x}^\beta)^{3/2}} = 0. \tag{4}$$

By $\varepsilon_{\alpha\beta}$ we denote here the Levi-Civita symbol $\begin{pmatrix} 0 & 1 \\ -1 & 0 \end{pmatrix}$. The Lagrange function for Eq. (4) has the form

$$\mathcal{L} = \frac{\varepsilon_{\alpha\beta} \dot{x}^\alpha \ddot{x}^\beta}{(\dot{x}_\beta \dot{x}^\beta)^{3/2}} - m\sqrt{\dot{x}_\beta \dot{x}^\beta}.$$

We now introduce the notation

$$\boldsymbol{u} = (u^\alpha) = (\dot{x}^\alpha) \tag{5}$$

and recall the definition of scalar product and the norm of two-vectors:

$$(\boldsymbol{a} \wedge \boldsymbol{b}) \cdot (\boldsymbol{w} \wedge \boldsymbol{z}) = (\boldsymbol{a} \cdot \boldsymbol{w})(\boldsymbol{b} \cdot \boldsymbol{z}) - (\boldsymbol{a} \cdot \boldsymbol{z})(\boldsymbol{b} \cdot \boldsymbol{w}),$$

$$\|\boldsymbol{a} \wedge \boldsymbol{b}\| = \sqrt{(\boldsymbol{a} \wedge \boldsymbol{b}) \cdot (\boldsymbol{a} \wedge \boldsymbol{b})},$$

or, in the coordinate form:

$$(\boldsymbol{a} \wedge \boldsymbol{b}) \cdot (\boldsymbol{w} \wedge \boldsymbol{z}) = 2(\boldsymbol{a} \wedge \boldsymbol{b})_{\alpha\beta} (\boldsymbol{w} \wedge \boldsymbol{z})^{\alpha\beta}$$

$$= a_{[\alpha} b_{\beta]} w^{[\alpha} z^{\beta]}.$$

In order to show that all integral curves of Eq. (4) have a constant curvature

$$k = \frac{\|u \wedge \dot{u}\|}{\|u\|^3}, \tag{6}$$

we first find the derivative of expression (6):

$$\dot{k} = \frac{(u \wedge \dot{u}) \cdot (u \wedge \ddot{u})}{\|u\|^3 \|u \wedge \dot{u}\|} - 3\frac{\|u \wedge \dot{u}\|(u \cdot \dot{u})}{\|u\|^5}. \tag{7}$$

We now solve Eq. (4) with respect to the variable $\ddot{u}$ by contracting it with the contravariant Levi-Civita symbol

$$e^{\alpha\beta} = \begin{pmatrix} 0 & 1 \\ -1 & 0 \end{pmatrix}:$$

$$\frac{\ddot{u}^\alpha}{\|u\|^3} = 3\frac{(u \cdot \dot{u})}{\|u\|^5}\dot{u}^\alpha + m\frac{(u \cdot u)e^{\beta\alpha}\dot{u}_\beta - (u \cdot \dot{u})e^{\beta\alpha}u_\beta}{\|u\|^3}. \tag{8}$$

Further, we substitute expression (8) in (7):

$$\dot{k} = 3\frac{(u \cdot \dot{u})(u \wedge \dot{u}) \cdot (u \wedge \dot{u})}{\|u\|^5 \|u \wedge \dot{u}\|} + \frac{m}{\|u\|^3}\Big[(u \wedge u)^2 e^{\beta\alpha}\dot{u}_\alpha\dot{u}_\beta$$

$$- (u \cdot u)(u \cdot \dot{u})(e^{\beta\alpha} + e^{\alpha\beta})u_\alpha\dot{u}_\beta + (u \cdot \dot{u})^2 e^{\beta\alpha}u_\beta u_\alpha\Big]$$

$$- 3\frac{\|u \wedge \dot{u}\|(u \cdot \dot{u})}{\|u\|^5} = 0.$$

The generalization to the Riemannian space was obtained in [4].

## 3. Three-Dimensional Space

In the three-dimensional space, one can seek a parameter-invariant variational equation for nonparametrized curves. The procedure used to deduce this equation was described in detail in [5]. It turns out that the problem of construction of a third-order parameter-invariant variational equation in a three-dimensional (pseudo)-Euclidean space invariant under the action of the (pseudo)-Euclidean group of transformations of this space is solved unambiguously. It is of interest to note that, in this case, it is not necessary to impose an *a priori* requirement of constancy of the Frénet curvature along the solution of the required equation. On the contrary, this property appears *automatically*. The parameter-invariant variational equation in variables (6) takes the form

$$\frac{\ddot{u} \times u}{\|u\|^3} - 3\frac{(\dot{u} \cdot u)}{\|u\|^5}\dot{u} \times u + m\frac{(u \cdot u)\dot{u} - (\dot{u} \cdot u)u}{\|u\|^3} = 0. \tag{9}$$

To show that the solutions of Eq. (9) have a constant Frénet curvature, we parametrize these solutions along the length of arc $s$ and denote the derivatives with respect to $s$ by the corresponding subscript:

$$u_s = u_s^\alpha = \dot{x}_s^\alpha = \frac{dx^\alpha}{ds}.$$

As a result, we obtain

$$u_s \cdot u_s = 1, \quad u_s \cdot \dot{u}_s = 0, \quad u_s \cdot \ddot{u}_s = -\dot{u}_s^2, \quad \dddot{u}_s \cdot u_s = -3\dot{u}_s \cdot \ddot{u}_s. \tag{10}$$

In these variables, Eq. (9) takes the form

$$\ddot{u}_s \times u_s + m\dot{u}_s = \mathbf{0}. \tag{11}$$

Moreover, the relation for curvature can be represented as

$$k = \|\dot{u}_s\|. \tag{12}$$

Further, we find the derivative of expression (12):

$$\frac{dk}{ds} = \frac{\dot{u}_s \cdot \ddot{u}_s}{\|\dot{u}_s\|}. \tag{13}$$

We now find the vector product of Eq. (11) by $\dot{u}_s$:

$$(\dot{u}_s \cdot \ddot{u}_s)u_s = \mathbf{0},$$

multiply this equation scalarly by $u_s$:

$$(\dot{u}_s \cdot \ddot{u}_s) = 0,$$

and substitute the result in relation (13), which gives $\dot{k} = 0$.

Further, we determine the value of torsion (second Frénet curvature)

$$\tau = \frac{u_s \cdot (\dot{u}_s \times \ddot{u}_s)}{k^2}.$$

Equation (11) immediately implies that

$$\tau = \frac{\dot{\boldsymbol{u}}_s \cdot (\ddot{\boldsymbol{u}}_s \times \boldsymbol{u}_s)}{k^2} = -m,$$

i.e., the integral curves are helical lines.

If $m = 0$, then these helical lines become plane and are called geodesic circles. For the Euclidean metric, these are ordinary circles of any radius. In the Lorentz metric, we have hyperbolas and, in the case of $(2 + 1)$-dimensional special theory of relativity, we speak about the world lines of relativistic uniformly accelerated motion.

## 4. Four-Dimensional Space-Time

The geometric facts about the variational problem in a four-dimensional space presented in what follows remain true both for the Euclidean and (pseudo)-Euclidean cases without serious changes. In view of the fact that the conclusions of the present work have interesting interpretations in the relativistic theory, we consider the case of Minkowski space-time with Lorentz metric

$$g_{\alpha\beta} = \text{diag}(1,-1,-1,-1). \tag{14}$$

In the previous works [2, 5], it was shown that Eq. (9) describes a spatially plane motion of a relativistic top with a constant spin of magnitude $\sigma$ (directed perpendicularly to the plane of motion) and a mass $m_0 = -m\sigma$. Thus, we consider the equations of motion of a free relativistic top in the Minkowski space of the special theory of relativity.

As shown in [2], the Mathisson–Papapetru equations used to describe the dynamics of a relativistic top with spin $\sigma^\alpha$, together with the so-called Mathisson–Pirani *additional condition*

$$\sigma^\alpha u_\alpha = 0, \tag{15}$$

take the following form in the Minkowski space:

$$\varepsilon_{\alpha\beta\gamma\mu}\ddot{u}^\beta u^\gamma \sigma^\mu - 3\frac{\dot{u}_\beta u^\beta}{u_\beta u^\beta}\varepsilon_{\alpha\beta\gamma\mu}\dot{u}^\beta u^\gamma \sigma^\mu - m_0\left(u_\beta u^\beta \dot{u}_\alpha - \dot{u}_\beta u^\beta u_\alpha\right) = 0, \tag{16}$$

where the four-vector of spin $\sigma^\alpha$ is now constant both in magnitude and in the direction.

**Proposition 1.** *Equation (16) is parameter-invariant.*

*Proof.* If we replace the parameter along the integral curve $\zeta \to \xi$, then the derivatives are transformed accordingly, $u_\zeta \to u_\xi$, by the following rule:

$$u_\zeta = \xi'_\zeta \, u_\xi,$$

$$\dot{u}_\zeta = \xi''_\zeta \, u_\xi + \xi'^2_\zeta \, \dot{u}_\xi,$$

$$\ddot{u}_\zeta = \xi'''_\zeta u_\xi + 3\xi''_\zeta \, \xi'_\zeta \, \dot{u}_\xi + \xi'^3_\zeta \, \ddot{u}_\xi.$$

We denote the left-hand side of Eq. (16) by $\mathcal{E}_\alpha(u^\beta, \dot{u}^\beta, \ddot{u}^\beta)$. By using the formulas presented above in Eq. (16), we obtain

$$\mathcal{E}_\alpha(u^\beta_\zeta, \dot{u}^\beta_\zeta, \ddot{u}^\beta_\zeta) = \xi'^4_\zeta \mathcal{E}_\alpha(u^\beta_\xi, \dot{u}^\beta_\xi, \ddot{u}^\beta_\xi).$$

Hence, we can choose the arc length $s$ as a parameter along the geodesic line. As a result, Eq. (16) takes a simpler form:

$$\varepsilon_{\alpha\beta\gamma\mu} \ddot{u}^\beta_s u^\gamma_s \sigma^\mu - m_0 \dot{u}_{s\alpha} = 0. \tag{17}$$

By using the contraction of Eq. (17) with the vector $\ddot{u}_s$, we can easily show that it has the curvature of the integral curve $k = \sqrt{\dot{u}_{s\alpha} \dot{u}^\alpha_s}$ as the first integral.

We now try to eliminate the (constant) vector variable $\sigma^\alpha$ by increasing the order of Eq. (17). However, we first contract Eq. (17) with a contravariant quantity $e^{\alpha\nu\rho\lambda} u_{s\rho}\sigma_\lambda$ with respect to the first index by using the well-known formula

$$\varepsilon_{\alpha\beta\gamma\mu} e^{\alpha\nu\rho\lambda} = \delta^{\nu\rho\lambda}_{\beta\gamma\mu}, \tag{18}$$

where $\delta^{\nu\rho\lambda}_{\beta\gamma\mu}$ is the generalized Kronecker symbol. To simplify the notation, we omit the subscript $s$ and obtain

$$\varepsilon_{\alpha\beta\gamma\mu} \ddot{u}^\beta u^\gamma \sigma^\mu e^{\alpha\nu\rho\lambda} u_\rho \sigma_\lambda = \delta^{\nu\rho\lambda}_{\beta\gamma\mu} \ddot{u}^\beta u^\gamma \sigma^\mu u_\rho \sigma_\lambda$$

$$= \ddot{u}^\nu u^\rho \sigma^\lambda u_\rho \sigma_\lambda + \ddot{u}^\lambda u^\nu \sigma^\rho u_\rho \sigma_\lambda + \ddot{u}^\rho u^\lambda \sigma^\nu u_\rho \sigma_\lambda$$

$$- \ddot{u}^\lambda u^\rho \sigma^\nu u_\rho \sigma_\lambda - \ddot{u}^\nu u^\lambda \sigma^\rho u_\rho \sigma_\lambda - \ddot{u}^\rho u^\nu \sigma^\lambda u_\rho \sigma_\lambda$$

$$= \sigma^\lambda \sigma_\lambda (\ddot{u}^\nu - \ddot{u}^\rho u_\rho u^\nu)$$

$$= \sigma^\lambda \sigma_\lambda (\ddot{u}^\nu + k^2 u^\nu).$$

Here, we have used the condition (10) of unitarity of the four-vector of velocity $u_s$ differentiated two times and relation (13). As a result, Eq. (17) takes the form

$$\|\sigma\|^2 (\ddot{u}_s^\nu + k^2 u_s^\nu) = m_0 e^{\alpha\nu\rho\lambda} \dot{u}_{s\alpha} u_{s\rho} \sigma_\lambda. \tag{19}$$

Further, we differentiate Eq. (19) and substitute (17) in the obtained equation instead of its right-hand side by taking into account the sign of determinant of the metric tensor, i.e.,

$$e_{\alpha\nu\rho\lambda} = \det(g_{\mu\beta}) \varepsilon_{\alpha\nu\rho\lambda} = -\varepsilon_{\alpha\nu\rho\lambda}, \tag{20}$$

$$\|\sigma\|^2 (\dddot{u}_s + k^2 \dot{u}_s) = m_0^2 \dot{u}_s. \tag{21}$$

The constant quantities in this equation can be also represented via the four-vector of linear momentum:

$$\mathcal{P}_\alpha = \frac{m_0}{\|u\|} u_\alpha + \frac{1}{\|u\|^3} \varepsilon_{\beta\gamma\lambda\alpha} \dot{u}^\beta u^\gamma \sigma^\lambda. \tag{22}$$

Based on relations (6), (18), and (20), we compute the length of the four-vector $\mathcal{P}$ as follows:

$$\mathcal{P}^2 = \mathcal{P}_\alpha \mathcal{P}^\alpha = -\varepsilon_{\beta\gamma\lambda\alpha} e^{\rho\nu\mu\alpha} \dot{u}^\beta u^\gamma \sigma^\lambda \dot{u}_\rho u_\nu \sigma_\mu = m_0^2 - k^2 \sigma^2. \tag{23}$$

It follows from relation (22) that

$$\frac{\mathcal{P} \cdot u}{\|u\|} = m_0. \tag{24}$$

We now make some assumptions. Thus, if we restrict ourselves to the branch of Eq. (16), where $\mathcal{P}^2 > 0$, then we conclude that the quantity

$$\omega^2 = -\frac{\mathcal{P}^2}{\sigma^2} = -\frac{m_0^2}{\sigma^2} + k^2 > 0 \tag{25}$$

because it follows from expression (15) that the four-vector of spin $\sigma$ (in the nontrivial case) must be spacelike, i.e., $\sigma^2 < 0$. Hence, equation (21) takes the form

$$\dddot{u}_s + \omega^2 \dot{u}_s = 0. \tag{26}$$

Equation (26) has the following corollary:

$$\ddot{u}_s \cdot u_s = 0. \tag{27}$$

In view of the condition of unitarity of the four-vector of velocity $u_s$ (10) differentiated three times, we arrive at one more corollary:

$$\frac{d}{ds} k^2 = \dot{u}_s \cdot \ddot{u}_s = 0. \tag{28}$$

However, this was supposed in Eq. (17).

Contracting Eq. (26) with the vector $\ddot{u}_s$, we get one more integral of motion, namely, $\ddot{u}_s^2$:

$$\frac{d}{ds} \ddot{u}_s^2 = 2\ddot{u}_s \cdot \dddot{u}_s = 0. \tag{29}$$

We now analyze the behavior of the torsion function $\tau$ on the solutions of Eq. (26). By using the well-known relation

$$\tau k^2 = \| u_s \wedge \dot{u}_s \wedge \ddot{u}_s \|, \tag{30}$$

in view of (10), we obtain

$$\tau^2 k^4 = \begin{vmatrix} u_s \cdot u_s & u_s \cdot \dot{u}_s & u_s \cdot \ddot{u}_s \\ \dot{u}_s \cdot u_s & \dot{u}_s \cdot \dot{u}_s & \dot{u}_s \cdot \ddot{u}_s \\ \ddot{u}_s \cdot u_s & \ddot{u}_s \cdot \dot{u}_s & \ddot{u}_s \cdot \ddot{u}_s \end{vmatrix}$$

$$= \begin{vmatrix} 1 & 0 & -k^2 \\ 0 & k^2 & k\dfrac{dk}{ds} \\ -k^2 & k\dfrac{dk}{ds} & \ddot{u}_s \cdot \ddot{u}_s \end{vmatrix}$$

$$= \ddot{u}_s^2 k^2 - k^6 - k^2 \left( \frac{dk}{ds} \right)^2. \tag{31}$$

By virtue of (28) and (29), this implies that $\tau$, just as $k$, is also an integral of motion, and the following relation is true:

$$\frac{\ddot{u}_s^2}{\dot{u}_s^2} = k^2 + \tau^2. \tag{32}$$

**Proposition 2.** *Equation (26) is equivalent to the following system (which simply means that we fix the integral of motion $\dfrac{\ddot{u}_s^2}{\dot{u}_s^2}$):*

$$\begin{cases} \ddot{\boldsymbol{u}}_s + \dfrac{\ddot{\boldsymbol{u}}_s^2}{\dot{\boldsymbol{u}}_s^2}\dot{\boldsymbol{u}}_s = 0, & \text{(I)} \\ \omega^2 = \dfrac{\ddot{\boldsymbol{u}}_s^2}{\dot{\boldsymbol{u}}_s^2}. & \text{(II)} \end{cases} \qquad (33)$$

***To prove*** this proposition, it is necessary to show that the expression $\dfrac{\ddot{\boldsymbol{u}}_s^2}{\dot{\boldsymbol{u}}_s^2}$ is indeed an integral of motion for Eq. (**I**) in system (33). We deduce relation (28) by repeating for Eq. (**I**) in (33) exactly the same operations as in the case of Eq. (26). Finally, we arrive at relation (29) by contracting Eq. (**I**) in (33) with the vector $\ddot{\boldsymbol{u}}_s$ and using the (already established) relation (28).

*4.1. Variational Formulation of the Mechanics of Constant Curvature in the Four-Dimensional Space-Time.* In this section, we show how to deduce a nonparametrized variational equation whose set of solutions contains the solutions of Eq. (26). It seems likely that the expression

$$\mathcal{R} = -\frac{1}{2}\left(\boldsymbol{u}\cdot\boldsymbol{u} - \frac{\dot{\boldsymbol{u}}\cdot\dot{\boldsymbol{u}}}{\omega^2}\right) \qquad (34)$$

proposed in [7] may play the role of Lagrange function for Eq. (26). {This was one of the first works on the so-called quiver motion ("Zitterbewegung") of quasiclassical spin, which is frequently cited even now [6]}. However this is not true. Equation (26) is written under the condition of validity of constraint (10). It can be easily obtained from the Lagrange function (34) if we forget that the subscript $s$ denotes a natural parameter, i.e., also forget about constraint (10) subjected to differentiation with respect to the natural parameter. In this case, the subscript $s$ can be omitted. However, the equation

$$\ddot{\boldsymbol{u}} + \omega^2\boldsymbol{u} = 0 \qquad (35)$$

obtained from the variational problem

$$\delta\int -\frac{1}{2}\left(\boldsymbol{u}\cdot\boldsymbol{u} - \frac{\dot{\boldsymbol{u}}\cdot\dot{\boldsymbol{u}}}{\omega^2}\right)d\zeta = 0, \qquad (36)$$

just as the variational problem (36) itself, is not nonparametrized, i.e., it gives no freedom of choosing the parameter $\zeta$ along its solutions and, hence, we have no right to treat $s$ as the independent variable in this equation. *Equation (26) is not equation (35) and, hence, cannot be obtained from the variational principle (36).* The variational problem with the Lagrange function (34) includes constraint (10). However, it is impossible to "transfer" this constraint directly to the equation of motion (35) and assume that the variable of integration is free as in problem (36). We now prove this fact.

Thus, we consider a variational problem with the Lagrange function (34) and a constraint

$$\psi = u_\alpha u^\alpha - 1. \qquad (37)$$

The variational equation takes the form

$$-\frac{d}{d\zeta}\left(\frac{\partial \mathcal{R}}{\partial \boldsymbol{u}}\right) + \frac{d^2}{d\zeta^2}\left(\frac{\partial \mathcal{R}}{\partial \dot{\boldsymbol{u}}}\right) - \frac{d}{d\zeta}\left(\lambda \frac{\partial \psi}{\partial \boldsymbol{u}}\right) = 0.$$

After multiplication by $\omega^2$, with regard for constraint (37), we obtain

$$\ddot{\boldsymbol{u}}_s + (1 - 2\lambda)\omega^2 \dot{\boldsymbol{u}}_s - 2\frac{d\lambda}{ds}\omega^2 \boldsymbol{u}_s = \boldsymbol{0}. \tag{38}$$

To determine $\lambda$ and $\dfrac{d\lambda}{ds}$, we contract Eq. (38) first with $\dot{\boldsymbol{u}}_s$ and then with $\boldsymbol{u}_s$ and apply relations (10):

$$1 - 2\lambda = -\frac{\ddot{\boldsymbol{u}}_s \cdot \dot{\boldsymbol{u}}_s}{\omega^2 \dot{\boldsymbol{u}}_s^2}, \tag{39}$$

$$\frac{d\lambda}{ds} = -\frac{3}{2}\frac{\ddot{\boldsymbol{u}}_s \cdot \dot{\boldsymbol{u}}_s}{\omega^2}. \tag{40}$$

Equation (38), with regard for (39) and (40), takes the form

$$\ddot{\boldsymbol{u}}_s - \frac{\ddot{\boldsymbol{u}}_s \cdot \dot{\boldsymbol{u}}_s}{\dot{\boldsymbol{u}}_s^2}\dot{\boldsymbol{u}}_s + 3(\ddot{\boldsymbol{u}}_s \cdot \dot{\boldsymbol{u}}_s)\boldsymbol{u}_s = \boldsymbol{0}. \tag{41}$$

Contracting with $\ddot{\boldsymbol{u}}_s$ and taking into account (10), we arrive at the following corollary of Eq. (41):

$$\ddot{\boldsymbol{u}}_s \cdot \ddot{\boldsymbol{u}}_s - \left(\frac{\ddot{\boldsymbol{u}}_s \cdot \dot{\boldsymbol{u}}_s}{\dot{\boldsymbol{u}}_s^2} + 3\dot{\boldsymbol{u}}_s^2\right)\ddot{\boldsymbol{u}}_s \cdot \dot{\boldsymbol{u}}_s \equiv 2\tau^2 k \frac{dk}{ds} + k^2 \tau \frac{d\tau}{ds} = 0. \tag{42}$$

The identity in this relation follows from the differential prolongations of relations (13) and (31), which yield the following respective expressions:

$$\ddot{\boldsymbol{u}}_s \cdot \boldsymbol{u}_s = k\frac{d^2 k}{ds^2} - k^2 \tau^2 - k^4, \tag{43}$$

$$\ddot{\boldsymbol{u}}_s \cdot \dot{\boldsymbol{u}}_s = k\tau^2 \frac{dk}{ds} + k^2 \tau \frac{d\tau}{ds} + 2k^3 \frac{dk}{ds} + \frac{dk}{ds}\frac{d^2 k}{ds^2}. \tag{44}$$

**Remark 1.** Relation (42) demonstrates that the curvature $k$ and torsion $\tau$ of the solutions of Eq. (41) can be constant only simultaneously.

Expressions (39) and (40) require compatibility in a sense that the second expression should be obtained as a result of differentiation of the first expression. In view of relation (44), this requirement imposes the following constraint:

$$\dddot{u}_s \cdot \dot{u}_s + 6(\ddot{u}_s \cdot \dot{u}_s)\dot{u}_s^2 - \frac{(\ddot{u}_s \cdot \dot{u}_s)(\ddot{u}_s \cdot \dot{u}_s)}{\dot{u}_s^2}$$

$$\equiv k\frac{d^3 k}{ds^3} - \frac{dk}{ds}\frac{d^2 k}{ds^2} + k^3 \frac{dk}{ds} - 2k\tau^2 \frac{dk}{ds} - 3k^2\tau\frac{d\tau}{ds} = 0. \quad (45)$$

The identity in this relation is obtained by using relation (44) and differential extension of relation (43).

**Proposition 3.** *The expression*

$$\frac{\tilde{A}}{2} = \frac{3}{2}\dot{u}^2 + \frac{\ddot{u}_s \cdot \dot{u}_s}{\dot{u}_s^2} \quad (46)$$

*is the first integral of the system of equations (41) and (45).*

*Proof.* We now differentiate (46) and apply (42) and (45).

*Remark 2.* Equation (41) is the *extension* of Eq. (**I**) of system (33) in a sense that, under the additional condition

$$\frac{dk}{ds} = 0$$

that does not contradict constraints (42) and (45), it is reduced to

$$\ddot{u}_s + \frac{\ddot{u}_s^2}{\dot{u}_s^2}\dot{u}_s = 0. \quad (47)$$

If we fix the value of curvature $k = k_0$, then, in view of (42), torsion $\tau$ takes a constant value $\tau_0$ and hence, according to relation (31), on the right-hand side of Eq. (47), we obtain

$$\frac{\ddot{u}_s^2}{\dot{u}_s^2} = k_0^2 + \tau_0^2$$

and, therefore, Eq. (41) is transformed into the equation of spin motion (26) with a frequency

$$\omega^2 = k_0^2 + \tau_0^2.$$

**4.2. Nonparametrized Variational Equation Containing the Equation of Spin Motion (26).** First, we consider the general form of variational equation for a nonparametrized variational problem. By setting $u^0 = 1$, every nonparametrized variational problem written in coordinates (5) can be rewritten in the following coordinates

$$t = x^0,$$

$$v^i = \frac{dx^i}{dt}, \quad v' = \frac{dv^i}{dt}, \quad v'' = \frac{dv'^i}{dt}, \quad v''' = \frac{dv''^i}{dt}, \quad i > 0. \tag{48}$$

The nonparametrized variational equation with left-hand side of the form (1) for the Lagrange function $\mathcal{L}(x^\alpha, u^\alpha, \dot{u}^\alpha)$ can be represented in variables (48) with the help of expression (2) according to the rule

$$\mathcal{E}_i = u^0 E_i, \tag{49}$$

by virtue of relation (3) and with the corresponding Lagrange function $L(t, x^i, v^i, v'^i)$ connected with the previous function by the formula

$$\mathcal{L}(x^\alpha, u^\alpha, \dot{u}^\alpha) = u^0 L(t, x^i, v^i, v'^i), \tag{50}$$

where the variables $t$, $v^i$ $v^i$, and $v'^i$ can be expressed via the variables $x^\alpha$, $u^\alpha$, and $\dot{u}^\alpha$.

Expression (2) autonomous with respect to time can be represented in the form

$$\mathbf{E} = \frac{d}{dt}\left(-\frac{\partial L}{\partial \mathbf{v}} + \frac{d}{dt}\frac{\partial L}{\partial \mathbf{v}'}\right). \tag{51}$$

We now study the problem of existence of expression (2) that generates expression (1) according to rule (49) so that, in the vector equation $\{\mathcal{E}_\alpha = 0\}$ independent of the parametrization of the integral curve, expression (1) obtained in passing to the natural parametrization $ds = \sqrt{1+\mathbf{v}^2}\, dt$ transforms into the first term of Eq. (26). In this case, we set $v_i = -v^i$ according to the expression of metric (15). For this purpose, we recalculate the fourth derivative

$$\ddot{u}^i_s = \frac{d}{ds}\ddot{u}^i_s = \frac{d}{ds}\frac{d}{dt}\ddot{u}^i_s \tag{52}$$

to the time variable according to which

$$\dot{t} = \frac{1}{\sqrt{1+\mathbf{v}^2}}, \quad \ddot{t} = -\frac{\mathbf{v}'\cdot\mathbf{v}}{(1+\mathbf{v}^2)^2}, \quad \dddot{t} = 4\frac{(\mathbf{v}'\cdot\mathbf{v})^2}{(1+\mathbf{v}^2)^{7/2}} - \frac{\mathbf{v}'^2 + \mathbf{v}''\cdot\mathbf{v}}{(1+\mathbf{v}^2)^{5/2}}.$$

We get

$$\dot{u}_s^i = \dot{t} v'^i, \quad \ddot{u}_s^i = \ddot{t} v'^i + \dot{t}^2 v''^i, \quad \dddot{u}_s^i = \dddot{t} v'^i + 3\dot{t}^2 \ddot{t} v''^i + \dot{t}^3 v'''^i; \quad (53)$$

in particular,

$$\ddot{u}_s^i = \left\{ 4 \frac{(\mathbf{v}' \cdot \mathbf{v})^2}{(1+\mathbf{v}^2)^{7/2}} - \frac{\mathbf{v}'^2 + (\mathbf{v}'' \cdot \mathbf{v})}{(1+\mathbf{v}^2)^{5/2}} \right\} v^i$$

$$- 3 \frac{(\mathbf{v}' \cdot \mathbf{v})}{(1+\mathbf{v}^2)^{5/2}} v'^i + \frac{v''^i}{(1+\mathbf{v}^2)^{3/2}}. \quad (54)$$

We now rewrite (54) in the form

$$\ddot{u}_{si} = \left\{ 4 \frac{(\mathbf{v}' \cdot \mathbf{v})^2}{(1+\mathbf{v}^2)^{7/2}} - \frac{\mathbf{v}'^2}{(1+\mathbf{v}^2)^{5/2}} - \frac{(\mathbf{v}'' \cdot \mathbf{v})}{(1+\mathbf{v}^2)^{5/2}} \right\} v_i$$

$$+ \frac{1}{2} \frac{d}{dt} \frac{\partial}{\partial v'^i} \frac{\mathbf{v}'^2}{(1+\mathbf{v}^2)^{3/2}}. \quad (55)$$

Further, we require the coincidence of (52) with (49), i.e., the validity of identity

$$\frac{d}{dt} \ddot{u}_s^i \equiv \frac{d}{dt} \left( -\frac{\partial L}{\partial v^i} + \frac{d}{dt} \frac{\partial L}{\partial v'^i} \right). \quad (56)$$

To this end, we first consider the identity

$$\ddot{u}_s^i \equiv -\frac{\partial L}{\partial v^i} + \frac{d}{dt} \frac{\partial L}{\partial v'^i}. \quad (57)$$

The comparison of (55) with (57) demonstrates that, for the last term in (55), it is reasonable to consider the Lagrange function

$$L_1 = \frac{1}{2} \frac{\mathbf{v}'^2}{(1+\mathbf{v}^2)^{3/2}}. \quad (58)$$

In identity (57), we take $L_1 + L_2$ instead of $L$, where $L_1$ is given by expression (58). Then the term $\dfrac{d}{dt} \dfrac{\partial L_1}{\partial v'^i}$ in (57) cancels with the last term in (55). As a result of collecting of the similar terms in (57), we can expand the operator

$$\frac{d}{dt} = v'^i \frac{\partial}{\partial v^i} + v''^i \frac{\partial}{\partial v'^i}. \quad (59)$$

This yields the following identity for $L_2$:

$$\left\{ 4\frac{(\mathbf{v}'\cdot\mathbf{v})^2}{(1+\mathbf{v}^2)^{7/2}} - \frac{5}{2}\frac{\mathbf{v}'^2}{(1+\mathbf{v}^2)^{5/2}} - \frac{(\mathbf{v}''\cdot\mathbf{v})}{(1+\mathbf{v}^2)^{5/2}} \right\} \mathbf{v}$$

$$\equiv -\frac{\partial L_2}{\partial \mathbf{v}} + v'^i \frac{\partial}{\partial v^i} \frac{\partial L_2}{\partial \mathbf{v}'} + v''^i \frac{\partial}{\partial v'^i} \frac{\partial L_2}{\partial \mathbf{v}'}.$$

The form of the coefficients at $v''^i$ suggests the following choice:

$$\frac{\partial L_2}{\partial \mathbf{v}'} = -\frac{(\mathbf{v}'\cdot\mathbf{v})}{(1+\mathbf{v}^2)^{5/2}},$$

which is integrated to give the following dependence of $L_2$ on $\mathbf{v}'$:

$$L_2 = -\frac{(\mathbf{v}'\cdot\mathbf{v})^2}{2(1+\mathbf{v}^2)^{5/2}}. \tag{60}$$

Continuing iterations, we substitute $L = L_1 + L_2 + L_3$ in the right-hand side of (57). As a result, the term

$$v''^i \frac{\partial}{\partial v'^i} \frac{\partial L_2}{\partial \mathbf{v}'}$$

in (57) cancels with the third term in expression (55). After collecting similar terms in (57), we substi-tute (57) in (56) and obtain the following identity for $L_3$:

$$-\frac{3}{2}\frac{d}{dt}\left\{\frac{\mathbf{v}'^2}{(1+\mathbf{v}^2)^{5/2}} - \frac{(\mathbf{v}'\cdot\mathbf{v})^2}{(1+\mathbf{v}^2)^{7/2}}\right\}\mathbf{v} \equiv -\frac{d}{dt}\left\{\frac{\partial L_3}{\partial \mathbf{v}} + \frac{d}{dt}\frac{\partial L_3}{\partial \mathbf{v}'}\right\}. \tag{61}$$

In this identity, the expression in braces on the left-hand side is nothing else but the square of the first Frénet curvature (12) with the factor $\frac{1}{\sqrt{1+v^2}}$. According to (28), the curvature $k$ must be constant, e.g., $k = k_0$, along the required variational equation. Identity (61) can be transformed as follows:

$$\frac{3}{2} k_0^2 \frac{d}{dt} \frac{\mathbf{v}}{\sqrt{1+\mathbf{v}^2}} \equiv \frac{d}{dt}\left\{\frac{\partial L_3}{\partial \mathbf{v}} + \frac{d}{dt}\frac{\partial L_3}{\partial \mathbf{v}'}\right\}. \tag{62}$$

Hence, in order to satisfy identity (62), it is sufficient to take the Lagrange function of free motion

$$L_3 = \frac{3}{2} k_0^2 \sqrt{1+\mathbf{v}^2}. \tag{63}$$

The Lagrange function for the second term in (26) [written with the time parameter according to (49) and (53)] is also the Lagrange function of free motion:

$$-\omega^2\sqrt{1+\mathbf{v}^2}. \tag{64}$$

Thus, it remains to combine (58), (60), (63), and (64) in order to get the Lagrange function for required variational equation

$$L = \frac{1}{2}(k^2 + A)\sqrt{1+\mathbf{v}^2},$$

where

$$A = 3k_0^2 - 2\omega^2.$$

In the homogeneous coordinates, according to relation (50) and definition (6), which has the same form in spaces of arbitrary dimension, we find

$$\mathcal{L} = \frac{1}{2}\left(\frac{u^2\dot{u}^2 - (u\cdot\dot{u})^2}{\|u\|^5} + A\|u\|\right). \tag{65}$$

### 4.3. Variational Extension of the Equation of Spin Motion.

**Proposition 4.** *The solutions of Eq. (26) are located among the extremals of the variational problem*

$$\delta\int\frac{1}{2}\left(\frac{u^2\dot{u}^2 - (u\cdot\dot{u})^2}{\|u\|^5} + A\|u\|\right)d\zeta = 0. \tag{66}$$

*Proof.* The variational equation for (66) is

$$\frac{d}{d\zeta}\left\{-2\frac{\ddot{u}}{\|u\|^3} + 6\frac{(u\cdot\dot{u})}{\|u\|^5}\dot{u}\right\} + \left(2\frac{(u\cdot\ddot{u})}{\|u\|^5} - \frac{\dot{u}^2}{\|u\|^5} - 5\frac{(u\cdot\dot{u})^2}{\|u\|^7} + \frac{A}{\|u\|}\right)u = \mathbf{0}. \tag{67}$$

Note that the variational problem (66) is not parametrized. Thus, in Eq. (67), we can pass to the natural parameter $s$ according to (10):

$$\ddot{u}_s + \left(\frac{3}{2}\dot{u}_s^2 - \frac{A}{2}\right)\dot{u}_s + 3(\dot{u}_s\cdot\ddot{u}_s)u_s = \mathbf{0}. \tag{68}$$

In view of relation (30), along the solutions of Eq. (68), we get

$$\frac{d\tau}{ds} = \frac{(u_s \wedge \dot{u}_s \wedge \ddot{u}_s) \cdot (u_s \wedge \dot{u}_s \wedge \dddot{u}_s)}{k^2 \|u_s \wedge \dot{u}_s \wedge \ddot{u}_s\|} - 2\frac{\|u_s \wedge \dot{u}_s \wedge \ddot{u}_s\|}{k^3}\dot{k}$$

$$= -2\frac{\tau}{k}\frac{dk}{ds}. \tag{69}$$

Note that, according to well-known formula

$$k^3 \tau^2 k_3 = \|u_s \wedge \dot{u}_s \wedge \ddot{u}_s \wedge \dddot{u}_s\|, \tag{70}$$

the third Frénet curvature $k_3$ is equal to zero along Eq. (68).

In order to choose solely the solutions of Eq. (26) in the collection of solutions of Eq. (68), it is necessary to impose an additional constraint

$$k^2 = \frac{2}{3}\omega^2 + \frac{A}{3}.$$

Relation (69) shows that, in addition to the constant curvature $k$, torsion $\tau$ also is an integral of motion. Hence, in this case, we speak about helical lines located in the three-dimensional space-time, i.e., about motion in a plane of the three-dimensional space.

In view of arguments presented above, it is reasonable to say that Eq. (67) is the *variational extension* of the equation of motion of quasiclassical spin (26) (as shown above, this equation is not variational).

### 4.4. Variational Extension of Eq. (47).

**Proposition 5.** *A system of equations*

$$\begin{cases} \dddot{u}_s + \left(\frac{3}{2}\dot{u}_s^2 - \frac{A}{2}\right)\dot{u}_s + 3(\dot{u}_s \cdot \ddot{u}_s)u_s = 0, & \text{(I)} \\ \frac{dk}{ds} = 0 & \text{(II)} \end{cases} \tag{71}$$

*is equivalent to the following system:*

$$\begin{cases} \dddot{u}_s + \frac{\ddot{u}_s^2}{\dot{u}_s^2}\dot{u}_s = 0, & \text{(I)} \\ k^2 - 2\tau^2 = A. & \text{(II)} \end{cases} \tag{72}$$

*Proof.* We first differentiate constraint (**II**) from system (71) and obtain

$$\ddot{u} \cdot \dot{u} = -\ddot{u}^2. \tag{73}$$

Further, we contract Eq. (**I**) in system (71) with the vector $\dot{u}$ and, in view of (73) and constraint (**II**) from system (71), obtain

$$\ddot{u}^2 = \left(\frac{3}{2}\dot{u}^2 - \frac{A}{2}\right)\dot{u}^2. \tag{74}$$

This expression is also substituted in Eq. (**I**) of system (71) in order to get Eq. (**I**) in system (72). In view of (31), Eq. (74), together with constraint (**II**) from system (71), is transformed into Eq. (**II**) in system (72).

Conversely, repeating the same transformations as earlier in the case of Eq. (32), we conclude that properties (28) and (29) are preserved for Eq. (**I**) in system (72). Hence, $k$ and $\tau$ are the first integrals, which immediately follows from (32). Therefore, constraint (**II**) in system (72) is consistent with the first equation of this system. Property (28) gives the second equation in system (71). Further, we substitute constraint (**II**) from system (72) in the form (74) in the first equation in (72) and arrive at the first two terms of the first equation in (71). The third term is identically equal to zero as a result of taking into account the already obtained constraint (**II**) from the same system (71).

*4.5. Variational Extension of the Lagrange Problem for Quasiclassical Spin.* On Eq. (68), the integral of motion (46) takes the value $\frac{A}{2}$. Constraint (45) is in this case unnecessary. Hence, the one-parameter family of Eqs. (67) can be regarded as the *variational extension* of the Lagrange problem (34), (37).

## CONCLUSIONS

We make an attempt to construct a variational principle for the paths of constant curvature in (pseudo)-Euclidean two-, three-, and four-dimensional spaces (inverse variational problem). Since the expression for curvature contains the second derivative, the required variational equations should be of the third order. In a two-dimensional space, the variational problem for which the equation of extremals contains the derivative of the third order (but not higher) can be formulated solely in the essentially parametric form. In this case, the inverse variational problem is solved unambiguously if we demand that the straight lines obtained as a particular case of solutions of the obtained equation should carry a natural parametrization by the parameter $s$ along them. In the three-dimensional space, there exists an absolutely unambiguous answer to the analyzed inverse variational problem according to which, despite the fact that the variational equation has (pseudo)-Euclidean symmetry, the Lagrange functions generating this equation [5] are not and cannot be invariants of the group of motions of the corresponding space [1]. In the four-dimensional space, it is impossible to construct a third-order variational equation with (pseudo)-Euclidean symmetry [5, III.3.7.5]. Generalizing the variational equation for the three-dimensional space to the four-dimensional space by introducing an additional dependence on a certain constant vector, we arrive at the relativistic variational equation of motion of the third order in the four-dimensional space-time, which is physically meaningful. Further, by applying the procedure of *elimination with the help of derivatives* of certain additional parameters, we obtain the variational equation in the four-dimensional (pseudo)-Euclidean space, which contains the equation of helical lines and generalizes the three-dimensional case. In this way, we get the variational generalization of differential equations by the method of differentiation with subse-quent "unfreezing" of the integrals of motion.

*Appendix.* **Illustration of Spin "Quiver."**

Equation (26) has a real solution

$$x^\mu = x_0^\mu + u_0^\mu s + a^\mu \cos(\omega s) + b^\mu \sin(\omega s). \tag{77}$$

The initial conditions are subjected to constraint (10).[1] In view of the independence of trigonometric binomials, this gives the following conditions:

$$a_\mu a^\mu = b_\mu b^\mu, \quad a_\mu b^\mu = 0, \quad a_\mu u_0^\mu = 0, \quad b_\mu u_0^\mu = 0. \tag{75}$$

Moreover, constraint (10) still adds the following additional condition:

$$u_{0\mu} u_0^\mu + a_\mu a^\mu \omega^2 = 1. \tag{76}$$

We know that, in space, the process of motion takes place in the plane $k_3 = 0$ and is described by relation (70). It is possible to consider the following values of the initial conditions:

$$a^1 = a^2 = b^1 = b^2 = \alpha, \quad a^0 = b^0 = \beta, \quad a^3 = b^3 = u_0^3 = 0. \tag{77}$$

We now get

$$a_\mu a^\mu = 2\alpha^2 - \beta^2.$$

However, it follows from the second relation in (75) that

$$a_\mu a^\mu = 2\alpha^2 - \beta^2 = 0. \tag{81}$$

Therefore, as a result of splitting into the space and time components, condition (76) turns into the following equality:

$$(u_0^0)^2 - (u_0^1)^2 - (u_0^2)^2 - (u_0^3)^2 = 1,$$

i.e., into the condition of unitarity for the initial velocity. However, in this case, the third relation in (75) implies that

$$a_\mu u_0^\mu = \alpha(u_0^1 + u_0^2) - \beta u_0^0 = 0. \tag{82}$$

---

[1] This condition was neglected in [7].

Thus, relation (81), together with (82), finally implies that initial conditions must satisfy the constraint

$$\frac{(u_0^1 + u_0^2)^2}{(u_0^0)^2} = 2.$$

Solution (77), split into the space and time components, takes the following final form:

$$\mathbf{r} = \mathbf{r}_0 + \mathbf{v}s + \alpha(\mathbf{i} + \mathbf{j})\big(\cos(\omega s) + \sin(\omega s)\big),$$

$$t = t_0 + \sqrt{2}\,(v_1 + v_2)\,s + \sqrt{2}\,\alpha\big(\cos(\omega s) + \sin(\omega s)\big),$$

where

$$\mathbf{r} = (x, y, z), \qquad \mathbf{v} = (v_1, v_2, v_3) \equiv (u_0^1, u_0^2, u_0^3), \qquad \mathbf{i} = (1, 0), \quad \mathbf{j} = (0, 1).$$

Note that the physical time $t$ cannot run proportionally to the proper time $s$ because, in this case, the amplitude $\alpha = a_1 = a_2$ would disappear.

*Dependence of the Character of Motion on the Ratio of Parameters.*

In Figs. 1–4, the dashed lines describe the process of motion without angular velocity $\omega = 0$.

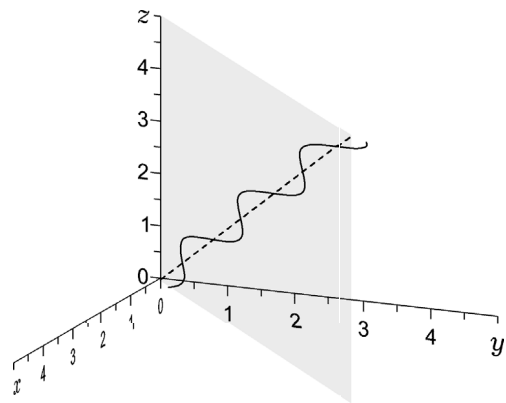

$v_1 = v_2 = v_3 = 1,$

$a_1 = a_2 = 1, \quad a_3 = 0, \quad \omega = 4.$

**Fig. 1**

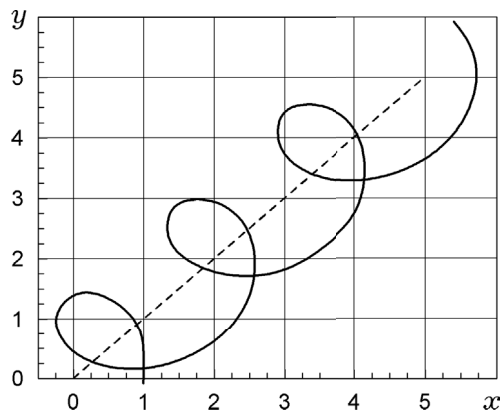

$v_1 = v_2 = v_3 = 1,$

$a_1 = a_2 = 0.3, \quad a_3 = 0, \quad \omega = 4.$

**Fig. 2**

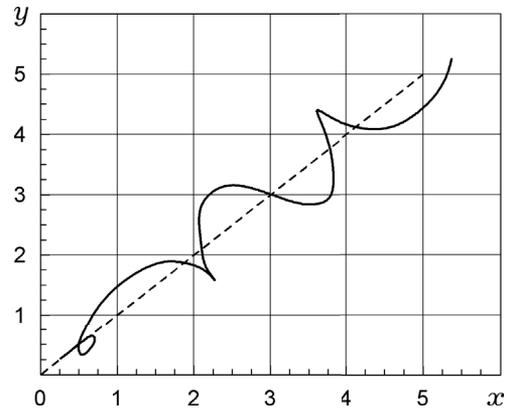

$$v_1 = v_2 = 10, \quad v_3 = 1,$$
$$a_1 = a_2 = 3, \quad a_3 = 0, \quad \omega = 1.52.$$

**Fig. 3**

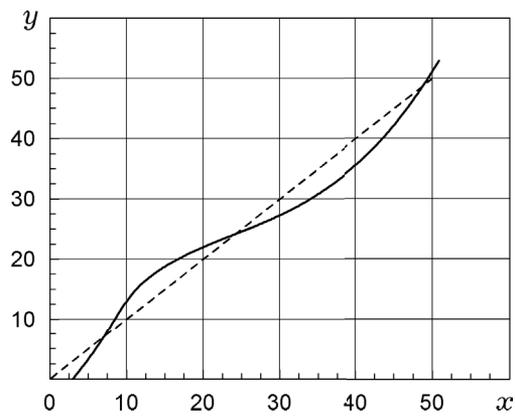

$$v_1 = v_2 = v_3 = 1,$$
$$a_1 = a_2 = 0.3, \quad a_3 = 0, \quad \omega = 5.$$

**Fig. 4**